\documentclass[10pt]{article}
\usepackage{cite}
\usepackage{mathrsfs}
\usepackage{amsfonts}
\usepackage{amsmath}
\usepackage{amsfonts,amssymb}
\usepackage{dsfont}
\usepackage{curves}
\usepackage{mathrsfs}
\usepackage{pifont}
\usepackage{amssymb}
\allowdisplaybreaks

\numberwithin{equation}{section}

\date{}

\def\BigRoman{\uppercase\expandafter{\romannumeral\number\count 255 }}
\def\Romannumeral{\afterassignment\BigRoman\count255=}

\begin{document}
\title{Spectral radius for the existence of $H_b$-factors in binding graphs
%\thanks{Supported by }
}
\author{\small Sufang Wang$^{1}$, Wei Zhang$^{2}$\footnote{Corresponding
author. E-mail address: wangsufangjust@163.com (S. Wang), zw\_wzu@163.com (W. Zhang).}\\
\small $1$.  School of Public Management, Jiangsu University of Science and Technology,\\
\small Zhenjiang, Jiangsu 212100, China\\
\small $2$.  School of Economics and Management, Wenzhou University of Technology,\\
\small Wenzhou, Zhejiang 325000, China
}

\maketitle
\begin{abstract}
\noindent The binding number, denoted by $\mbox{bind}(G)$, of a graph $G$ is defined as the minimum value of $\frac{|N_G(X)|}{|X|}$ taken over any non-empty subset $X$ of $V(G)$ with $N_G(X)\neq V(G)$.
A graph $G$ is said to be $r$-binding if $\mbox{bind}(G)\geq r$. The adjacency matrix of a graph $G$ is denoted by $A(G)$. The largest eigenvalue of $A(G)$ is called the spectral radius of $G$. An
$H_b$-factor of a graph $G$ is defined as a spanning subgraph $F$ of $G$ such that for any $v\in V(G)$, $d_F(v)$ belongs to the set $\{1,3,5,\ldots,b-1,b\}$, where $b$ is an even integer with $b\geq2$.
This note establishes a sufficient condition to guarantee that a connected $\frac{1}{b-1}$-binding graph $G$ of even order contains an $H_b$-factor based on the spectral radius.
\\
\begin{flushleft}

(2020) Mathematics Subject Classification: 05C50, 05C70, 90B99
\end{flushleft}
\end{abstract}

\section{Introduction}

Graphs discussed here are simple and undirected. Let $G$ be a graph with vertex set $V(G)=\{v_1,v_2,\ldots,v_n\}$ and edge set $E(G)$. The order of $G$ is denoted by $|V(G)|=n$. For $v_i\in V(G)$,
the degree of $v_i$, denoted by $d_G(v_i)$, is the number of edges adjacent to $v_i$ in $G$. Let $o(G)$ denote the number of odd components in $G$. For any $S\subseteq V(G)$, let $G[S]$ denote the
subgraph of $G$ induced by $S$, and let $G-S$ denote the subgraph obtained from $G$ by deleting the vertices in $S$ together with their incident edges. The complete graph of order $n$ is denoted by
$K_n$. For two vertex disjoint graphs $G_1$ and $G_2$, we use $G_1\cup G_2$ and $G_1\vee G_2$ to denote the disjoint union and the join of $G_1$ and $G_2$, respectively. The adjacency matrix of a
graph $G$ is defined as the $(0,1)$-matrix $A(G)=(a_{ij})_{n\times n}$, where $a_{ij}=1$ if and only if $v_iv_j\in E(G)$. The largest eigenvalue of $A(G)$, denoted by $\rho(G)$, is called the
spectral radius of $G$. The reader can be referred to \cite{PB,ELW,Wc,Ws1,Za,ZZZL,WZL} for some properties of the spectral radius in $G$.

For any vertex $v$ of $G$, the neighborhood of $v$ in $G$ is denoted by $N_G(v)$. For any vertex subset $X$ of $G$, let $N_G(X)=\bigcup\limits_{v\in X}N_G(v)$. Woodall \cite{Wt} introduced the concept
of the binding number. The binding number, denoted by $\mbox{bind}(G)$, of $G$ is defined as the minimum value of $\frac{|N_G(X)|}{|X|}$ taken over any non-empty subset $X$ of $V(G)$ with $N_G(X)\neq V(G)$.
A graph $G$ is said to be $r$-binding if $\mbox{bind}(G)\geq r$.

Let $b$ be an even integer with $b\geq2$. An odd $[1,b-1]$-factor of $G$ is a spanning subgraph $F$ of $G$ such that $d_F(v)$ is odd and $1\leq d_F(v)\leq b-1$ for any $v\in V(G)$. In particular, an
odd $[1,1]$-factor is also called a perfect matching. An $H_b$-factor of $G$ is defined as a spanning subgraph $F$ of $G$ such that for any $v\in V(G)$, $d_F(v)$ belongs to the set $\{1,3,5,\ldots,b-1,b\}$.

O \cite{O} verified a lower bound for the spectral radius to guarantee the existence of a perfect matching in a connected graph. Fan and Lin \cite{FL} determined a spectral condition to guarantee that
a 1-binding graph contains a perfect matching. Fan, Lin and Lu \cite{FLL} provided a spectral condition for the existence of an odd $[1,b-1]$-factor in a graph with minimum degree. Zhou \cite{Zs1}
proposed two sufficient conditions for a graph to possess an odd $[1,b-1]$-factor with given property. Lu and Wang \cite{LW} showed some results on the existence of an $H_b$-factor in a connected graph.
Fan and Liu \cite{FL1} established a lower bound based on the spectral radius of a graph $G$ to ensure that $G$ contains an $H_b$-factor. Zhou \cite{Zs2} presented some sufficient conditions with respect
to the size and the signless Laplacian spectral radius for a graph with minimum degree to possess an $H_b$-factor. Wu \cite{Wd} provided an upper bound in terms of the distance spectral radius of a graph
$G$ to guarantee the existence of an $H_b$-factor in $G$. Some other sufficient conditions on the existence of graph factors were provided by Ananchuen, Caccetta and Ananchuen \cite{ACA}, Egawa, Kano and
Yan \cite{EKY}, Wu \cite{Wu}, Pan and Zhou \cite{PZ}, Zhou \cite{Zs3,Zhou}, Zhou, Bian and Sun \cite{ZBS}, Zhou, Bian and Wu \cite{ZBW}, Zhou, Zhang and Sun \cite{ZZS}, O \cite{Oe}, Kim and O \cite{KO}.

Motivated by \cite{LW} directly, it is natural and interesting to put forward a sufficient condition based on the spectral radius to ensure that a graph has an $H_b$-factor. Our main result is given in
the following.

\medskip

\noindent{\textbf{Theorem 1.1.}} Let $b$ be an even integer with $b\geq2$, and let $G$ be a connected $\frac{1}{b-1}$-binding graph of even order $n\geq b+6$. If
$$
\rho(G)\geq\rho(K_1\vee(K_{n-b-3}\cup K_3\cup(b-1)K_1)),
$$
then $G$ contains an $H_b$-factor unless $G=K_1\vee(K_{n-b-3}\cup K_3\cup(b-1)K_1)$.

\medskip

Based on Theorem 1.1, the following corollary holds.

\medskip

\noindent{\textbf{Corollary 1.2.}} Let $G$ be a connected 1-binding graph of even order $n\geq8$. If
$$
\rho(G)\geq\rho(K_1\vee(K_{n-5}\cup K_3\cup K_1)),
$$
then $G$ contains a $[1,2]$-factor unless $G=K_1\vee(K_{n-5}\cup K_3\cup K_1)$.

\medskip

\section{Some preliminaries}

Lu and Wang \cite{LW} showed a sufficient condition for a graph with even order to possess an $H_b$-factor.

\medskip

\noindent{\textbf{Lemma 2.1}} (Lu and Wang \cite{LW}). Let $b$ be an even integer with $b\geq2$, and let $G$ be a connected graph of even order. If
$$
o(G-S)\leq b|S|
$$
for any nonempty subset $S$ of $V(G)$, then $G$ contains an $H_b$-factor.

\medskip

\noindent{\textbf{Lemma 2.2}} (Li and Feng \cite{LF}). Let $H$ be a subgraph of a connected graph $G$. Then
$$
\rho(G)\geq\rho(H),
$$
with equality following if and only if $G=H$.

\medskip

\noindent{\textbf{Lemma 2.3}} (Fan and Lin \cite{FL}). Let $s+\sum\limits_{i=1}^{t}n_i=n$ with $s\geq1$. If $n_t\geq n_{t-1}\geq\cdots\geq n_1\geq1$, $n_{t-1}\geq3$ and $n_t\leq n-s-t-1$, then
$$
\rho(K_s\vee(K_{n_t}\cup K_{n_{t-1}}\cup\cdots\cup K_{n_1}))\leq\rho(K_s\vee(K_{n-s-t-1}\cup K_3\cup(t-2)K_1)),
$$
with equality holding if and only if $(n_t,n_{t-1},n_{t-2},\ldots,n_1)=(n-s-t-1,3,1,\ldots,1)$.

\medskip

Let $M$ be an $n\times n$ real matrix and $V=\{1,2,\ldots,n\}$. Given a partition $\pi: V=V_1\cup V_2\cup\cdots\cup V_r$, the matrix $M$ can be correspondingly partitioned as
\begin{align*}
M=\left(
  \begin{array}{cccc}
    M_{11} & M_{12} & \cdots & M_{1r}\\
    M_{21} & M_{22} & \cdots & M_{2r}\\
    \vdots & \vdots & \ddots & \vdots\\
    M_{r1} & M_{r2} & \cdots & M_{rr}\\
  \end{array}
\right),
\end{align*}
where $M_{ij}$ denotes the submatrix (block) of $M$ formed by rows in $V_i$ and columns in $V_j$. Let $b_{ij}$ be the average row sum of $M_{ij}$. The quotient matrix of $M$ with respect to $\pi$ is defined
as the matrix $B_{\pi}=(b_{ij})_{r\times r}$. The partition $\pi$ is called equitable if every block $M_{ij}$ of $M$ admits constant row sum $b_{ij}$.

\medskip

\noindent{\textbf{Lemma 2.4}} (You, Yang, So and Xi \cite{YYSX}). Let $M$ denote a real $n\times n$ matrix with an equitable partition $\pi$, and let
$M_{\pi}$ be the corresponding quotient matrix. Then the eigenvalues of $M_{\pi}$ are also eigenvalues of $M$. Furthermore, if $M$ is nonnegative
and irreducible, then the largest eigenvalues of $M$ and $M_{\pi}$ are equal.

\medskip

\section{The proof of Theorem 1.1}

\medskip

\noindent{\it Proof of Theorem 1.1.} Suppose that $G$ contains no $H_b$-factor. According to Lemma 2.1, there exists a nonempty subset $S\subseteq V(G)$ satisfying $o(G-S)\geq b|S|+1$. Set $|S|=s$ and $o(G-S)=q$.
Then $q\geq bs+1$. The $q$ odd components in $G-S$ are denoted by $O_1,O_2,\ldots,O_q$. We write $|O_i|=n_i$ for $1\leq i\leq q$. Without loss of generality, let $n_q\geq n_{q-1}\geq\cdots\geq n_1$. We first
prove the following claim.

\medskip

\noindent{\bf Claim 1.} $n_{bs}\geq3$.

\noindent{\it Proof.} Assume that $n_{bs}=1$. Obviously, we possess $n_i=1$ for $1\leq i\leq bs$ due to $1\leq n_1\leq n_2\leq\cdots\leq n_{bs}$. Let $X=V(O_1\cup O_2\cup\cdots\cup O_{bs})$. Then we easily see
$N_G(X)\subseteq S$ and $|X|=bs$. Thus, we have
$$
\frac{|N_G(X)|}{|X|}\leq\frac{|S|}{|X|}=\frac{s}{bs}=\frac{1}{b}<\frac{1}{b-1},
$$
which contradicts that $G$ is $\frac{1}{b-1}$-binding. Hence, we deduce $n_{bs}\geq3$. Claim 1 is proved. \hfill $\Box$

\medskip

Obviously, $n_i\geq3$ for $q\geq i\geq bs$ due to Claim 1, $q\geq bs+1$ and $n_q\geq n_{q-1}\geq\cdots\geq n_1$. We easily see that $G$ is a spanning subgraph of
$G_1=K_s\vee(K_{n'_{bs+1}}\cup K_{n_{bs}}\cup K_{n_{bs-1}}\cup\cdots\cup K_{n_1})$, where $n'_{bs+1}\geq n_{bs}\geq\cdots\geq n_1\geq1$, $n_{bs}\geq3$ and $n'_{bs+1}=n-s-\sum\limits_{i=1}^{bs}n_i$. In terms of
Lemma 2.2, we infer
\begin{align}\label{eq:3.1}
\rho(G)\leq\rho(G_1),
\end{align}
with equality following if and only if $G=G_1$.

Let $G_2=K_s\vee(K_{n-(b+1)s-2}\cup K_3\cup(bs-1)K_1)$, where $n\geq(b+1)s+5$. Using Lemma 2.3, we conclude
\begin{align}\label{eq:3.2}
\rho(G_1)\leq\rho(G_2),
\end{align}
with equality occurring if and only if $(n'_{bs+1},n_{bs},n_{bs-1},\ldots,n_1)=(n-(b+1)s-2,3,1,\ldots,1)$. The following proof will be divided into two cases.

\noindent{\bf Case 1.} $s=1$.

Recall that $G_1=K_s\vee(K_{n'_{bs+1}}\cup K_{n_{bs}}\cup K_{n_{bs-1}}\cup\cdots\cup K_{n_1})$. Applying $s=1$ and Lemma 2.3, we obtain
$$
\rho(G_1)\leq\rho(K_1\vee(K_{n-b-3}\cup K_3\cup(b-1)K_1)),
$$
where the equality follows if and only if $G_1=K_1\vee(K_{n-b-3}\cup K_3\cup(b-1)K_1)$. Combining this with \eqref{eq:3.1}, we get
$$
\rho(G)\leq\rho(K_1\vee(K_{n-b-3}\cup K_3\cup(b-1)K_1)),
$$
where the equality holds if and only if $G=K_1\vee(K_{n-b-3}\cup K_3\cup(b-1)K_1)$, a contradiction.

\noindent{\bf Case 2.} $s\geq2$.

Recall that $n\geq(b+1)s+5$ and $G_2=K_s\vee(K_{n-(b+1)s-2}\cup K_3\cup(bs-1)K_1)$. Based on the partition $V(G_2)=V(K_s)\cup V(K_{n-(b+1)s-2})\cup V(K_3)\cup V((bs-1)K_1)$, the quotient matrix of $A(G_2)$ can
be written as
\begin{align*}
B_2=\left(
  \begin{array}{cccc}
  s-1 & n-(b+1)s-2 & 3 & bs-1\\
  s & n-(b+1)s-3 & 0 & 0\\
  s & 0 & 2 & 0\\
  s & 0 & 0 & 0\\
  \end{array}
\right),
\end{align*}
and the characteristic polynomial of $B_2$ equals
\begin{align*}
\psi_{B_2}(x)=&x^{4}+(-n+bs+2)x^{3}+(n-bs^{2}-(b+2)s-5)x^{2}\\
&+((bs^{2}+2s+2)n-(b^{2}+b)s^{3}-(3b+2)s^{2}-(2b+8)s-6)x\\
&-(2bs^{2}-2s)n+(2b^{2}+2b)s^{3}+(4b-2)s^{2}-6s.
\end{align*}
In light of Lemma 2.4 and the equitable partition $V(G_2)=V(K_s)\cup V(K_{n-(b+1)s-2})\cup V(K_3)\cup V((bs-1)K_1)$, the largest root of $\psi_{B_2}(x)=0$ is equal to $\rho(G_2)$.

Let $G_*=K_1\vee(K_{n-b-3}\cup K_3\cup(b-1)K_1)$. Then the quotient matrix, corresponding to the partition $V(G_*)=V(K_1)\cup V(K_{n-b-3})\cup V(K_3)\cup V((b-1)K_1)$, of $A(G_*)$ is
\begin{align*}
B_*=\left(
  \begin{array}{cccc}
  0 & n-b-3 & 3 & b-1\\
  1 & n-b-4 & 0 & 0\\
  1 & 0 & 2 & 0\\
  1 & 0 & 0 & 0\\
  \end{array}
\right).
\end{align*}
The characteristic polynomial of $B_*$ equals
\begin{align*}
\psi_{B_*}(x)=&x^{4}+(-n+b+2)x^{3}+(n-2b-7)x^{2}\\
&+((b+4)n-b^{2}-6b-16)x-(2b-2)n+2b^{2}+6b-8.
\end{align*}
By Lemma 2.4 and the equitable partition $V(G_*)=V(K_1)\cup V(K_{n-b-3})\cup V(K_3)\cup V((b-1)K_1)$, the largest root of $\psi_{B_*}(x)=0$ equals $\rho(G_*)$. Notice that
$G_*=K_1\vee(K_{n-b-3}\cup K_3\cup(b-1)K_1)$ contains $K_{n-b-2}$ as its proper subgraph. By Lemma 2.2, we infer
\begin{align}\label{eq:3.3}
\rho(G_*)>\rho(K_{n-b-2})=n-b-3.
\end{align}

Based on a simple calculation, we get
\begin{align}\label{eq:3.4}
\psi_{B_*}(x)-\psi_{B_2}(x)=(s-1)f(x),
\end{align}
where $f(x)=-bx^{3}+(bs+2b+2)x^{2}+(-(bs+b+2)n+(b^{2}+b)s^{2}+(b^{2}+4b+2)s+b^{2}+6b+10)x+(2bs+2b-2)n-(2b^{2}+2b)s^{2}-(2b^{2}+6b-2)s-2b^{2}-6b+8$. Then we easily see
$$
f'(x)=-3bx^{2}+2(bs+2b+2)x-(bs+b+2)n+(b^{2}+b)s^{2}+(b^{2}+4b+2)s+b^{2}+6b+10.
$$
When $x\geq n-b-3$, we deduce
\begin{align*}
f'(x)=&-3bx^{2}+2(bs+2b+2)x-(bs+b+2)n+(b^{2}+b)s^{2}+(b^{2}+4b+2)s+b^{2}+6b+10\\
\leq&-3b(n-b-3)^{2}+2(bs+2b+2)(n-b-3)-(bs+b+2)n\\
&+(b^{2}+b)s^{2}+(b^{2}+4b+2)s+b^{2}+6b+10\\
=&-3bn^{2}+(bs+6b^{2}+21b+2)n+(b^{2}+b)s^{2}-(b^{2}+2b-2)s-3b^{3}-21b^{2}-37b-2\\
\leq&-3b((b+1)s+5)^{2}+(bs+6b^{2}+21b+2)((b+1)s+5)\\
&+(b^{2}+b)s^{2}-(b^{2}+2b-2)s-3b^{3}-21b^{2}-37b-2\\
=&-(3b^{3}+4b^{2}+b)s^{2}+(6b^{3}-4b^{2}-4b+4)s-3b^{3}+9b^{2}-7b+8\\
\leq&-4(3b^{3}+4b^{2}+b)+2(6b^{3}-4b^{2}-4b+4)-3b^{3}+9b^{2}-7b+8\\
=&-3b^{3}-15b^{2}-19b+16\\
<&0
\end{align*}
due to $b\geq2$, $s\geq2$ and $n\geq(b+1)s+5$. This leads to
\begin{align}\label{eq:3.5}
f(x)\leq f(n-b-3)
\end{align}
for $x\geq n-b-3$. Then we compute
\begin{align}\label{eq:3.6}
f(n-b-3)=&-b(n-b-3)^{3}+(bs+2b+2)(n-b-3)^{2}\nonumber\\
&+(-(bs+b+2)n+(b^{2}+b)s^{2}+(b^{2}+4b+2)s+b^{2}+6b+10)(n-b-3)\nonumber\\
&+(2bs+2b-2)n-(2b^{2}+2b)s^{2}-(2b^{2}+6b-2)s-2b^{2}-6b+8\nonumber\\
=&-bn^{3}+(3b^{2}+10b)n^{2}+((b^{2}+b)s^{2}+(3b+2)s-3b^{3}-20b^{2}-30b+2)n\nonumber\\
&-(b^{3}+6b^{2}+5b)s^{2}-(3b^{2}+11b+4)s+b^{4}+10b^{3}+30b^{2}+23b-4\nonumber\\
\triangleq&\varphi(n).
\end{align}
In light of $b\geq2$, $s\geq2$ and $n\geq(b+1)s+5$, we possess
\begin{align*}
\varphi'(n)=&-3bn^{2}+2(3b^{2}+10b)n+(b^{2}+b)s^{2}+(3b+2)s-3b^{3}-20b^{2}-30b+2\\
\leq&-3b((b+1)s+5)^{2}+2(3b^{2}+10b)((b+1)s+5)\\
&+(b^{2}+b)s^{2}+(3b+2)s-3b^{3}-20b^{2}-30b+2\\
=&-(3b^{3}+5b^{2}+2b)s^{2}+(6b^{3}-4b^{2}-7b+2)s-3b^{3}+10b^{2}-5b+2\\
\leq&-4(3b^{3}+5b^{2}+2b)+2(6b^{3}-4b^{2}-7b+2)-3b^{3}+10b^{2}-5b+2\\
=&-3b^{3}-18b^{2}-27b+6\\
<&0.
\end{align*}
This yields that $\varphi(n)$ is decreasing for $n\geq(b+1)s+5$. Thus, we have
\begin{align}\label{eq:3.7}
\varphi(n)\leq&\varphi((b+1)s+5)\nonumber\\
=&-(b^{4}+2b^{3}+b^{2})s^{3}+(3b^{4}-5b^{2}+2)s^{2}\nonumber\\
&-(3b^{4}-7b^{3}-2b^{2}-b-8)s+b^{4}-5b^{3}+5b^{2}-2b+6\nonumber\\
\triangleq&g(s).
\end{align}
For $s\geq2$, it follows from $b\geq2$ that
\begin{align*}
g'(s)=&-3(b^{4}+2b^{3}+b^{2})s^{2}+2(3b^{4}-5b^{2}+2)s-(3b^{4}-7b^{3}-2b^{2}-b-8)\\
\leq&-12(b^{4}+2b^{3}+b^{2})+4(3b^{4}-5b^{2}+2)-(3b^{4}-7b^{3}-2b^{2}-b-8)\\
=&-3b^{4}-17b^{3}-30b^{2}+b+16\\
<&0.
\end{align*}
This leads to $g(s)\leq g(2)$ for $s\geq2$. Combining this with $b\geq2$, we obtain
$$
g(s)\leq g(2)=-b^{4}-7b^{3}-19b^{2}+30<0.
$$
Combining this with \eqref{eq:3.4}, \eqref{eq:3.5}, \eqref{eq:3.6}, \eqref{eq:3.7} and $s\geq2$, we deduce
$$
\psi_{B_*}(x)-\psi_{B_2}(x)=(s-1)f(x)\leq(s-1)f(n-b-3)=(s-1)\varphi(n)\leq(s-1)g(s)<0.
$$
Thus, we conclude $\psi_{B_*}(x)<\psi_{B_2}(x)$ for $x\geq n-b-3$. Together with \eqref{eq:3.1}, \eqref{eq:3.2} and \eqref{eq:3.3}, we have
$$
\rho(G)\leq\rho(G_1)\leq\rho(G_2)<\rho(G_*)=\rho(K_1\vee(K_{n-b-3}\cup K_3\cup(b-1)K_1)),
$$
which is a contradiction to $\rho(G)\geq\rho(K_1\vee(K_{n-b-3}\cup K_3\cup(b-1)K_1))$. Theorem 1.1 is proved. \hfill $\Box$

\medskip

\section*{Declaration of competing interest}

\medskip

The authors declare that they have no known competing financial interests or personal relationships that could have appeared to influence the work
reported in this paper.

\section*{Data availability}

\medskip

No data was used for the research described in the article.

\medskip

%\section*{Acknowledgments}


\begin{thebibliography}{9999}


\bibitem {PB} S. Pirzada, Z. Bhat, Brualdi-Solheid problem on the minimum spectral radius of graphs with given matching number, Discrete Math. 349(6) (2026) 115031.

\bibitem {ELW} M. Ellingham, L. Lu, Z. Wang, Maximum spectral radius of outerplanar 3-uniform hypergraphs, J. Graph Theory 100(4) (2022) 671--685.

\bibitem {Wc} J. Wu, Characterizing spanning trees via the size or the spectral radius of graphs, Aequationes Math. 98(6) (2024) 1441--1455.

\bibitem {Ws1} J. Wu, Some results on the $k$-strong parity property in a graph, Comput. Appl. Math. 45(4) (2026) 138.

\bibitem {Za} S. Zhou, A result on spanning trees with bounded total excess, Discrete Appl. Math. 388 (2026) 130--135.

\bibitem {ZZZL} S. Zhou, Y. Zhang, T. Zhang, H. Liu, Toughness and $A_{\alpha}$-spectral radius in graphs, Filomat 40(5) (2026) 1883--1892.

\bibitem {WZL} J. Wu, S. Zhou, H. Liu, A spectral condition for spanning trees with restricted degrees in bipartite graphs, Proc. Rom. Acad. Ser. A Math. Phys. Tech. Sci. Inf. Sci. 27(1) (2026) 19--24.

\bibitem {Wt} D. Woodall, The binding number of a graph and its Anderson number, J. Combin. Theory Ser. B 15 (1973) 225--255.

\bibitem {O} S. O, Spectral radius and matchings in graphs, Linear Algebra Appl. 614 (2021) 316--324.

\bibitem {FL} D. Fan, H. Lin, Binding number, $k$-factor and spectral radius of graphs, Electron. J. Combin. 31(1) (2024) \#P1.30.

\bibitem {FLL} D. Fan, H. Lin, H. Lu, Spectral radius and $[a,b]$-factors in graphs, Discrete Math. 345 (2022) 112892.

\bibitem {Zs1} S. Zhou, Sufficient conditions for a graph with minimum degree to be $k$-critical with respect to an odd $[1,b]$-factor, Discrete Math. 349(11) (2026) 115251.

\bibitem {LW} H. Lu, D. Wang, On Cui-Kano's characterization problem on graph factors, J. Graph Theory 74 (2013) 335--343.

\bibitem {FL1} A. Fan, R. Liu, Spectral radius and factors in graphs with given minimum degree, Appl. Math. Comput. 529 (2026) 130151.

\bibitem {Zs2} S. Zhou, Sufficient conditions for a special factor in a graph with minimum degree, https://arxiv.org/pdf/2606.08502.

\bibitem {Wd} J. Wu, Distance spectral radius and $H_b$-factors in graphs, https://arxiv.org/pdf/2606.20692.

\bibitem {ACA} N. Ananchuen, L. Caccetta, W. Ananchuen, A characterization of maximal non-$k$-factor-critical graphs, Discrete Math. 307 (2007) 108--114.

\bibitem {EKY} Y. Egawa, M. Kano, Z. Yan, $(1,f)$-Factors of graphs with odd property, Graphs Comb. 32 (2016) 103--110.

\bibitem {Wu} J. Wu, Sufficient conditions for a graph with minimum degree to have a component factor, Proc. Rom. Acad. Ser. A Math. Phys. Tech. Sci. Inf. Sci. 27(1) (2026) 3--10.

\bibitem {PZ} Q. Pan, S. Zhou, Sufficient conditions for isolated tough graphs to have path-factors, AIMS Math. 11(5) (2026) 13371--13383.

\bibitem {Zs3} S. Zhou, Some spectral conditions for star-factors in bipartite graphs, Discrete Appl. Math. 369 (2025) 124--130.

\bibitem {Zhou} S. Zhou, Spanning subgraphs and spectral radius in graphs, Aequationes Math. 100(1) (2026) 1.

\bibitem {ZBS} S. Zhou, Q. Bian, Z. Sun, Spectral conditions for path-factors in isolated tough graphs, Discrete Appl. Math. 385 (2026) 228–236.

\bibitem {ZBW} S. Zhou, Q. Bian, J. Wu, Sufficient conditions for even factors in graphs, Discrete Appl. Math. 386 (2026) 365--372.

\bibitem {ZZS} S. Zhou, Y. Zhang, Z. Sun, The $A_{\alpha}$-spectral radius for path-factors in graphs, Discrete Math. 347(5) (2024) 113940.

\bibitem {Oe} S. O, Eigenvalues and $[a,b]$-factors in regular graphs, J. Graph Theory 100 (2022) 458--469.

\bibitem {KO} D. Kim, S. O, Eigenvalues and parity factors in graphs with given minimum degree, Discrete Math. 346 (2023) 113290.

\bibitem {LF} Q. Li, K. Feng, On the largest eigenvalue of a graph, Acta Math Appl Sin Chinese Ser 2 (1979) 167--175.

\bibitem {FL} D. Fan, H. Lin, Binding number, $k$-factor and spectral radius of graphs, Electron. J. Combin. 31(1) (2024) \#P1.30.

\bibitem {YYSX} L. You, M. Yang, W. So, W. Xi, On the spectrum of an equitable quotient matrix and its application, Linear Algebra Appl. 577 (2019) 21--40.


\end{thebibliography}
\end{document}